# The Exact Spectral Derivative Discretization Finite Difference (ESDDFD) Method for Wave Models: A Wave View of the Universe Through Natural Fractional and Fractal Derivative Representations (or *View Lens Shops for The Exponential Wave Universe*)

D.P. Clemence-Mkhope
Department of Mathematics and Statistics
North Carolina A&T State University
Greensboro, NC 27411, USA
e-mail: clemence@ncat.edu

**Abstract**
A wave view of the universe is proposed in which each natural phenomenon is equipped with its own unique natural viewing lens. A self-sameness modeling principle and its systematic application in Fourier-Laplace transform space is proposed as a novel, universal discrete modeling paradigm for advection-diffusion-reaction equations (ADREs) across non-integer derivatives, time scales, and wave spectral signatures. Its implementation is a novel exact spectral derivative discretization finite difference method (ESDDFD), a way for crafting wave viewing lenses by obtaining discrete wave models from ADRE models. The template for building these lenses come in the form of natural derivative representations obtained from the wave signature probability distribution function and its harmonic oscillation in FL transform space; use of the ESDDFD method in the discrete numerical modeling of wave equations requires no a-priori theory of any mathematical derivative. A major mathematical consequence of this viewpoint is that all notions of the mathematical integer or non-integer derivatives have representation as limits of such natural derivative representations; this and other consequences are discussed and a discretization of a simple integer derivative diffusion-reaction equation is presented to illustrate the method. The resulting view lenses, in the form of ESDDFD models, work well in detecting both local and non-local Debye or Kohlrausch-Williams-Watts exponential patterns; only Brownian motion and sub-diffusion are discussed in the present article.

**Key Words:** Debye exponential patterns, Kohlrausch-Williams-Watts stretched exponential patterns, Brownian motion, sub-diffusion, fractional or fractal derivatives, nonstandard finite difference (NSFD) methods, advection-diffusion-reaction wave equations

# 1 Introduction

This is the introductory chapter of a thesis undertaken by the author in response to the following question posed, in 2008, in private correspondence by Professor Mickens (a copy of the letter is in the Thesis) about the order of three numerical schemes for the decay equation:

> **Question:** *What is the "order of accuracy" for the following three schemes for the ODE*
>
> (*) $\quad \dfrac{dx}{dt} = -\lambda x$
>
> (1) $\dfrac{x_{k+1} - x_k}{h} = -\lambda x_k$
>
> (2) $\dfrac{x_{k+1} - x_k}{h} = -\lambda x_{k+1}$
>
> (3) $\dfrac{x_{k+1} - x_k}{\phi(h)} = -\lambda x_{k+1}, \ \phi(h) = \dfrac{1}{\lambda}\left(e^{\lambda h} - 1\right)$
>
> *I would like to get your reasoning and the details of these calculations.*

These questions were in a way a culmination of a conversation started during the 2002 Mathematics Joint Meetings in New Orleans, wherein for a few days our ponderings are reflected in scribbled notes about exponential nonstandard (NSFD) denominators and classic equations, but mostly the Airy equation (a copy of the notes is in the Thesis, with a brief account of what led to that conversation).

Mickens was at the time advancing his NSFD methodology [Mick1] and explaining how it could be used to construct unconditionally convergent numerical models, including (2) and (3) above. My interest at the time was in exploring if numerical models of certain Sturm-Liouville systems (specifically Schrödinger and Dirac waves; see, e.g., [Fyn] for a concise primer on classical (Newtonian) and relativistic (Einsteinian) wave viewpoints as well as the non-relativistic (Schrödinger/Heisenberg/Born) and relativistic (Dirac) explanation of the quantum universe) could be used to 'see' half-bound states, whose tracking had been theorized (by my PhD Thesis advisor and his collaborators for the Schrödinger case on the half-line or whole-line [HS], the Dirac case on the half-line [4], and by myself on the whole line for the Dirac case [HKS]) as being possible through the behavior of the Titchmarsh-Weyl -$m(\lambda)$ function, a spectral descriptor ([C&L] has a gentle introduction to the $m(\lambda)$, including the underlying theories of ordinary differential equations (ODEs) and complex variables). While I agreed that his method seemed to have merit, I believed that it was still lacking something for realistic modeling of wave signatures: while it guaranteed that the resulting discrete model would retain essential features of the system under study, its formulation did not include any spectral information – Mickens' view of my pursuit of spectral theory is in the Thesis.

This Thesis is my effort to answer the question posed in that letter and to explain my reasoning in the best way I know. It is therefore a proposition severely limited by limited grasp of the laws that govern the interactions of phenomena in our universe and the mathematical nuances that purport to explain them. My background as a mathematical physicist limits me to an explanation based on a view of science as a study of wave phenomena through the viewing windows of differential equations.

My Thesis is that within each window, as represented by the various differential equation modeling paradigms (integer, fractional, fractal, stochastic, infinitesimal, etc.) on the ordinary or other time scales, there exists a natural spectral viewing lens for each phenomenon and a universal template for constructing that lens from its observed spectral (wave, wavelet, or string) signature pattern. It is embodied in the fundamental spectral behavior and its representation is invariant across all DE modeling paradigms. To capture wave behavior at spectral level, exact representation of this fundamental spectral behavior must form the basis for DE wave model discretization; hence the Exact Spectral Derivative Discretization Finite Difference (ESDDFD) method.

While the original motivation was to write an explanation to justify an answer to the questions posed in the Mickens letter, further motivation has been added along the way, especially by the dismissive attitude about the importance of studying equation (*) above as well as about the relevance of some fractional derivatives, particularly those of conformable type (see [Khl] and extensions described in, e.g., [Musr] [Vld], [Almd, [Imb], [H&K]), as descriptors of information other than that already described by integer derivatives. The following quotes reflect those sentiments and my opinions about them:

> '' Local fractional derivatives of differentiable functions are Integer-order derivatives or zero. … Therefore we can state that the conformable fractional derivatives do not give anything new in the spaces of differentiable functions and are not fractional derivatives of non-integer order." (see, e.g., [TEV1])

> "As one can observe from (the formula $_0^C T_t^\alpha (f)(t) = t^{1-\alpha} \frac{d}{dt} f(t)$) the so-called conformable fractional derivative (CFD) is nothing else than the integer-order derivative multiplied by a power function. There is no need to explicitly consider this derivative (which is actually not a fractional-order operator) since any differential problem with the CFD can be rewritten as a standard integer-order differential equation." [Private communication, 2021; see also [TEV2]

As is argued in [DPC-M1], these are (partly true but) mistaken points of view that may be a result of limited understanding of how fractional models arise in applications such as wave modeling; the conclusion in [DPC-M4] of the consequence of such viewpoints was the following:

> Since such a term as $t^\alpha$ in modeling may arise (see [M&K], for brief wave modeling accounts and references therein or [B&B] (see also [TEV3])for detailed exposition) from modeling observed behavior of that form, formulation of problems in terms of the first derivative leaves a whole class of phenomena unstudied: it is like viewing phenomena for $0 < \alpha < 1$ through lenses developed for $\alpha = 1$; no matter how good those lenses, some observations will be a little off!!

In the context of the proposed universal wave view, this is akin to viewing sub-diffusion processes through lenses crafted for Brownian motion (see, e.g., [M&K] and Chapters 40-42 of [Fyn]): they only see Debye exponential wave patterns and are blind to all Kohlrausch-Williams-Watts (KWW) stretched exponential wave patterns in the sub-diffusion regime.

It is this latter motivation that has resulted in an exposition emphasizing wave modeling and its link to non-integer derivatives (NIDs). The narration reference to 'wave lens shops' and the connection of various mathematical and scientific concepts is also an artifact of the history of this thesis: my late mother, in 2004, tasked me to explain, in a way she and her grandchildren could understand, what calculus had to do with knowing that the morning star (Venus) is a planet or the milky way is a galaxy (how that came about is in the thesis). This is therefore my thesis on the scientific theory of spectral discrete wave modeling to my friend and mentor, Ronald E. Mickens, and an explanation of the real-world applications of that theory to all my children (all Maritela's children and their children) and all my students (and their students) as a tribute to, and in memory of my mother, Nessie M. Mkhope.

The rest of this chapter is organized as follows. In the next section, the context for the proposed universal wave view is established, along with the basic principles for modeling in that universe. Differential equations (DE) modeling of diffusion-advection-reaction (ADR) processes is briefly described in Section 3. The discrete modeling of sub-diffusion DE models of ADR processes is briefly discussed in Section 4, wherein descriptions are given of the traditional, Mickens NSFD, and new ESDDFD methods. Section 5 presents a discretization of a linear, integer derivative example to illustrate the new method. To conclude, Section 6 presents twelve immediate and foreseen consequences of the ESDDFD method, which constitute an outline of the remaining chapters of the thesis.

**2 The Proposed Universal Wave View and its Spectral Modelling Principles**

**2.1 A Wave View of the Universe**

Information about various phenomena is received in observations in the form of waves – such as light, heat, sound, water; one goal of mathematical modeling is to understand the nature of those waves through the study of their signature patterns. This supports scientific progress in that it creates the ability to formulate recreation of aspects of the phenomena through observations of those signature patterns. It is akin to standing at the origin in Fig 1, looking at all the incoming lines as signature waves and, without knowing their origin, trying to recreate the formulas that produced them, hence identifying their origin.

One way to think about the science industry is therefore as the wave industry; its pursuits comprise the conceptualization, design, manufacture, maintenance, monitoring, and continuous upgrades of various instruments for wave signature detection. The modeling of wave signatures using mathematics and various laws from laboratory and other observations is just a small sub-industry of this conglomerate, whose products and by-products are turned into instruments to improve human life, and the use of differential equations to do so is yet a small, lens design component of that wave sub-industry.

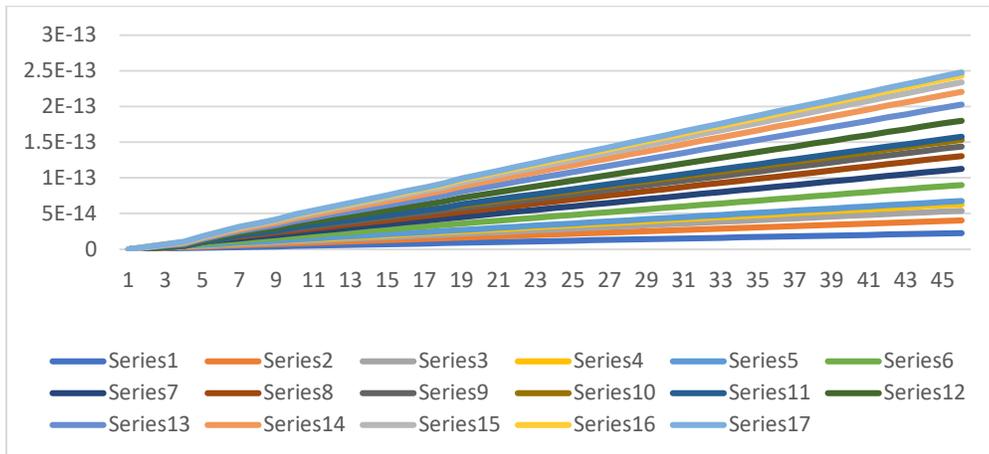

Figure 1. The Problem: How do we distinguish these same-looking wave-lines?

Through observations, signatures of wide classes of phenomena have been classified as obeying two power law patterns:

(I) The Debye exponential pattern, whose asymptotic wave signature behavior at the origin obeys a linear law:

$$W(t) \sim Ct, \tag{2.1.1}$$

for some constant C, or, more generally

$$W(t) \sim (\psi(t)) \tag{2.1.2}$$

where $\psi(t)$ is a well-behaved linear function determined by the process under study.

(II) The Kohlrausch-Williams-Watts (KWW) stretched exponential pattern, whose asymptotic wave signature behavior at the origin obeys a power law:

$$W(t) \sim Ct^\alpha, \tag{2.1.3}$$

for some constants $C, \alpha$ (including fractions) or, more generally,

$$y(t) \sim (\psi(t, \alpha)) \tag{2.1.4}$$

where $\psi(t, \alpha)$ is a well-behaved power function determined by the process under study.

Each of the behaviors (I) and (II) above, which according to [M&K], are "ubiquitous to a diverse number of systems", can further be classified into non-periodic and periodic wave signatures, where each type can be local or non-local; these behaviors are discussed further in Section 4.3

### 2.2 Proposed Universal Wave View Modeling Principles

It is the proposition of this thesis that the key to identifying each one of the signatures in (I) and (II) out of the ocean of signatures (an almost infinite number in mathematical thinking; I say 'almost infinite' because in reality they are not infinite – just because we cannot count them, it does not mean that a cup, or ocean, of (pure?) water contains an infinite number of $H_2O$ molecules – there may be Avogadro's number raised to whatever power you wish, but there are still a finite number of them!) is contained in two of what seem to be foundational principles for viewing our universe through the classical and newer theories of waves, each wrapped in one of two of the most basic differential equations, just in a more rudimentary form than often taught in differential equation courses.

The key divergence of the present proposition from traditional practice is the view that differential equations (DEs) as models are mathematical approximations of the original modeling templates, which are difference equations with spectral denominators, not Taylor truncations. DE models can be obtained as a limit of models derived from various scientific principles as difference equations describing process evolution using rates of change; while the traditional view treats these rates of change as approximations of tangent lines, the spectral view treats them as exact difference quotients in spectral space and is embodied in the following foundational principle for the proposed paradigm.

**The Exact Spectral Difference Discrete (ESDD) Wave Modeling Principle**

The propagation of all wave signatures is governed at Fourier-Laplace spectral level (that is, at what is currently accepted as the most fundamental level) by two self-sameness principles, which may be expressed as un-encumbered (I) wave relaxation or growth and (II) harmonic wave motion and represented by difference equations SSP(F) below. All derivative representations must be obtained from the exact form of these principles at that level to capture spectral details of the model.

The ESDFD method is a recognition that since, as stated in [Chetal], "The relaxation oscillation equation is the master equation of relaxation and oscillation processes.", and since this mastery is at spectral level, the best way to discretize derivatives is to do so at spectral level using this master equation.

**SSP(F) Exact *Self-sameness Modeling Principles as Rules for Rates of Change***

***SSP(FI) The first self-sameness principle***; the exact decay/growth equation: $\quad \frac{\Delta y(t)}{\mu_1(\Delta t, \lambda, \alpha)} = \pm \lambda y(t)$

***SSP(FII) The second self-sameness principle***; the exact HO equation: $\quad \frac{\Delta(\Delta y(t))}{(\mu_2(\Delta t, \omega, \alpha))^2} \pm \omega^2 y(t) = 0$

In SSP(F) above, $\Delta y(t) \equiv y(t) - y(x)$ is the change of the quantity $y$ between the points $t$ and $x$ as viewed from the point $t$ while $\mu_1(\Delta t, \lambda, \alpha)$ and $\mu_2(\Delta t, \omega, \alpha)$ denote the unique measures for each viewing each wave signature (viewing lens, as it were) that are crafted from the key; $\lambda$ and $\omega$ contain the specific wave information for each wave signature in that ocean. An interesting consequence of our approach is that $\Delta y(t)$ with respect to the correct measure is independent of the point $x$.

An immediate corollary to the (ESDD) Wave Modeling Principle is that (in 1-D) Debye and KWW time patterns ($\sim Ct^\alpha, 0 < \alpha \leq 1$) originate in Fourier transform [Z&C] space, while Debye and Kohlrausch-Williams-Watts space patterns ($\sim Cx^\beta, 0 < \beta \leq 1$) originate in Laplace transform [Z&C] space.

The equations embodying the principles are given in SSP(L) below in the form common in textbooks, where $\lambda$ may be positive or negative to represent a growth or decay rate and $\omega$ is a wave frequency.

**SSP(L) Limit *Self-sameness Modeling* Principles**

***SSP(LI) The first self-sameness principle***; the decay/growth equation: $\quad y' = \pm \lambda y$

***SSP(LII) The second self-sameness principle***; the harmonic oscillator equation: $\quad y'' \pm \omega^2 y = 0$

Physically, the mathematical limit approximation SSP(L) of the original model SSP(F) appears to represent a template for viewing lenses where all the signatures appear the same. It may be a great tool for identifying common characteristics of waves classes, but it is a very poor one for distinguishing between the individuals in those classes.

While it is beyond the scope of the present discussion, it is conjectured that some version of the same keys should be applicable to other views, such as those whose premise is that information signatures are received in the form of wavelets or strings. Instead of basing the wave tracking denominators on Fourier-Laplace spectra, it is reasonable to conjecture that similar use of the wavelet and string theory transforms will lead to parallel applications based on their versions of SSP(I) and SSP(II). Put another way, there should be the same template for the wavelet transform [W&W] as for the Fourier transform, with its versions of the Euler (exponential) and Mittag-Leffler (ML) functions ([Grnf], [Gmz]) as the basic building blocks for derivative representations for wavelet models; this is discussed further in Chapter 4.

Since our approach applies equally to Debye and KWW patterns (that is, Brownian motion as well as sub-diffusion), the generalized self-sameness principles may be stated as in SSP(G) below to reflect the fractional (see [Pdlb],[Mai]) or fractal (see, e.g., [Hlf], [Ch] and references therein) dimension); it is therefore assumed henceforth that $\alpha \in (0, 1]$:

**SSP(G). Generalized *Self-sameness Modeling* Principles**

***SSP(GI) The first self-sameness principle***; the exact decay/growth equation: $\quad \dfrac{\Delta_t^\alpha y(t)}{\mu_1(\Delta t, \lambda, \alpha)} = \pm \lambda y(t)$

***SSP(GII) The second self-sameness principle***; the exact HO equation: $\quad \dfrac{\Delta_t^\alpha \left( \Delta_t^\alpha y(t) \right)}{\left( \mu_2(\Delta t, \omega, \alpha) \right)^2} \pm \omega^2 y(t) = 0$

Since the exact propagator (that is, the solution, of SSP(GII) can be expressed in terms of that of SSP(GI), in what follows, SSP(GI) will often be referred to as the fundamental self-sameness principle.

Perhaps the best way to think of 'fractal space' (as opposed to fractional) is as regular space with holes, in the sense illustrated in [He]: as if a rectangle were composed of many very small rectangles separated by extremely small gaps. The smallest gap is the limiting constant for the step size in the He derivative (the He view treats fractal space just as if it were regular space, with a step size restriction to prevent ending up in one of the separating gaps). In the Hölder/ Hausdorff/Chen views (see, e.g.,

[Hlf],[Ch]), fractal space is transformed through a change of variables and the transformed space looks like regular space, whence derivative constructions in Hölder/Hausdorff space looks the same as that of regular space.

In the window-lens point of view, my thesis is that to date, to my best knowledge, the best universal lens template for viewing diverse classes of wave signatures through their spectral modes has been designed by Euler and Mittag-Leffler: Euler's template is just a local version of the non-local ML one. This template translates un-altered in form through various transformations, diffusion regimes space dimensions, wave origins (deterministic or stochastic), and spectral measures. While only the Brownian motion and sub-diffusion regimes are presented here, extension to super-diffusion is conjectured based on known extensibility of non-integer derivatives from $0 < \alpha \leq 1$ to $\alpha \geq 1$, which includes ballistic diffusion. An interesting feature of this lens template is that the resulting lens for each wave signature has its own incoming-outgoing (represented as decaying and growing) wave switch, explicitly given by two discrete representations of SSP(GI); SSP(GII), which models 'pure' periodic waves, does not differentiate between incoming and outgoing waves.

An interesting consequence of our approach, which needs no a-priori mathematical definition of a derivative to discretize differential equations, is that just like their integer cousin, NIDs arise as rather sanitized limit cases of exact difference quotients defined in terms of the solution of Eqn (*). Moreover, the conformable derivative is the link between local and non-local dynamics, leading to a possible expression of the ML function in terms of the Euler function; the ESDDFD method may therefore lead to further understanding about the switch between Debye and KWW patterns.

## 3 DE Modeling of Advection-Diffusion-Reaction (ADR) Brownian Motion and Sub-diffusion Processes

In all the cases mentioned above, the modeling process that results in the governing DEs does not require an a-priori definition of the derivative and involves no limit process. A rate of change is usually modeled and then its 'small denominator limit' is labeled as the derivative, which happens to coincide with the slope of the tangent line for $\alpha = 1$. For ease of exposition, the following presentation is limited to waves propagating in one space dimension; two- and three-dimensional extensions are given in [DPC-M4], which is Chapter 4 of the thesis.

DE modeling of ADR processes (see, e.g., [Fyn], [W&W], or [M&K] for an advanced primer) can be generally described by the following general partial differential equation (PDE), to which we will refer as the ADR equation (ADRE), with the functions $W, F(x, W), R(r, W)$ as well as the parameters $D$ and $r$ defined according to application:

$$\frac{\partial W}{\partial t} = \left[F(x,W)\frac{\partial}{\partial x} + D\frac{\partial^2}{\partial x^2} + R(r,W)\right] W(x,t). \tag{3.1}$$

The terms $\frac{\partial W}{\partial t}, F(x,W)\frac{\partial W}{\partial x}, D\frac{\partial^2 W}{\partial x^2}$, and $WR(r,W)$ in general represent, respectively, the time evolution of the system, the advection, diffusion, and reaction effects of the interaction on the process. In the ADRE (Eqn (3.1)), the $D$ being positive, negative, or zero represent, respectively decaying, growing, or 'standing' waves, whose models are mathematically studied as parabolic, elliptic, or hyperbolic problems: our method applies to all of them and is extensible to NIDEs and higher space dimensions.

The ADRE is usually derived (or explained, see [M&K], for example) as a first and second order truncations of a Taylor series for the first and second derivatives in either variable, which is based on the first order Taylor term being coincident with the slope of a tangent line obtained as a limit of secant lines. While the ADRE is normally derived as limit case, a version of it that is free of the limit process is the following difference equation, wherein all derivatives in the ADRE are replaced by their difference quotients:

$$\frac{\Delta W}{\mu(\Delta t,\lambda,\alpha)} = \left[F(x,W)\frac{\Delta}{\mu_1(\Delta x,\lambda,\alpha)} \pm D\frac{\Delta^2}{(\mu_2(\Delta x,\lambda,\alpha))^2} + R(r,W)\right]W(x,t), \qquad (3.2)$$

where the symbol $\Delta$ means the difference as defined in SSP(G) and not the Laplacian. In applications, the functions $F(x,W)$ and $R(r,W)$ satisfy some specified small-time asymptotic behavior. In what follows, we will write the ADRE generally as

$$H\left(W,\frac{\partial W}{\partial t},\frac{\partial W}{\partial x},\frac{\partial^2 W}{\partial x^2}\right) = 0, \text{ or } L\left(W,\frac{\partial W}{\partial t},\frac{\partial W}{\partial x},\frac{\partial^2 W}{\partial x^2}\right) = 0. \qquad (3.3)$$

(The $H-L$ notation for me is a reminder that these methods apply whether you arrive at Eq (3.1) by Hamiltonian or by Lagrangian. In practice, various initial and boundary values are attached, with conditions to assure mathematical well-posedness of the resulting initial value or boundary value problem (IVP or BVP, or both), that is, the existence of a unique solution for each set of specified conditions. Various types of these conditions are classified (such as the familiar Cauchy, Dirichlet, Neumann, or Robin conditions [Z&C]) according to the nature of the phenomenon under study, including the geometry within which it occurs. These IVP/BVP constructions result in what mathematicians then study as dynamical systems, which evolve given initial data according to the model: to be considered realistic, it is expected that the trajectory of this dynamic evolution will re-produce some known behavior of the phenomenon. Good models should be invariant under various transformations (such as those of Lorentz (see, e.g., [Fyn], Chapter 15, for a primer), Hausdorff (see, e.g., [K&G1][Ch]), and others) of coordinate systems and are expected to detect both classical as well as quantum effects, whether the view frame be relativistic or non-relativistic. In Brownian motion terms, this simply means that for each starting point and time, the trajectory of the model should re-create the random walk the original system purports to model. If 'small' changes in the initial point lead to only 'small' changes in the trajectory, the model is labeled (see, e.g. [JDL])as 'stable with respect to initial conditions' (meaning that it is good), otherwise it is considered chaotic, generating apparently random patterns.

From our view windows and lenses perspective, all this essential translates to the following statement about models: While different effects might be detected from the different viewpoints, an object viewed with Dirac (see, e.g., [Fyn], pg. 52-10) or Schrodinger/Heisenberg/Born lenses (who, according to [Fyn], pg. 37-1, "obtained a consistent description of the behavior of matter on a small scale") through the Einsteinian view window should 'be the same' as if viewed through the Newtonian window; different viewpoints of the same object should not see different objects, just different aspects of the same one – because it is the same object. Modeling propagation of an object's wave signature in exact form at spectral level uniquely distinguishes it amongst all others: the ESDDFD accomplishes exact discrete representation of the relaxation and oscillation rules at the Fourier-Laplace spectral level, thereby modeling distinctive signatures at spectral level.

While the focus of this thesis is mainly on (3.4), the time-fractional generalization of model (3.3),

$$H\left(W, D_t^\alpha, \frac{\partial W}{\partial x}, \frac{\partial^2 W}{\partial x^2}\right) = 0, 0 < \alpha \leq 1, \tag{3.4}$$

the principles and constructions presented are easily extensible to the more general time and space fractional or fractal systems,

$$H\left(W, D_t^\alpha W, \nabla_x^\beta W, \nabla_{xx}^{2\beta} W\right) = 0, 0 < \alpha, \beta \leq 1. \tag{3.5}$$

Analogous to the ADRE (3.1) and (3.2), the system (3.4) is a limit case of the following finite difference system (3.5), which our ESDDFD approach considers as the starting point, thereby omitting the derivative limit process from the construction of discrete models:

$$H\left(W, \Delta_t^\alpha W, \Delta_x^\beta W, \Delta_x^\beta\left(\Delta_x^\beta W\right)\right) = 0, 0 < \alpha, \beta \leq 1, \tag{3.6}$$

In this more general case, the replacement of $D_t^\alpha, \nabla_x^\beta$ and $\nabla_{xx}^{2\beta}$ are their exact representations, respectively from SSP(GI) and SSP(GII) in Fourier, Laplace, or Fourier-Laplace transform space.

The ESDDFD representation replacements for Eqn (3.6) neither depend on nor require any knowledge about the mathematical notion or theory of derivative or integral beyond Fourier and Laplace transform tables. This diverges from all currently published approaches for the discretization of non-integer derivatives (NIDs) by traditional Euler or NSFD methods. The method discussed here applies to NID differential equations (NIDEs) for a wide variety of NIDs: the non-local classic fractional derivatives of Caputo and Riemann-Liouville (see, e.g., [Mai],[Pdlb]) and their recent extensions ([CF], [A&B]), the classical fractal derivatives of Hölder/Hausdorff type (see, e.g.,[Hlf]) and their recent extensions ([K&G],[Ch], [He])), the local conformable fractional derivative (CFD) ([Khl]) and its extensions and generalizations (e.g., [Musr], [Vld], [KCi],[Z&LC][Almd], and a good summary in [Imb]), including, as is shown in [DPC-M2], their stochastic (see., e.g [GJ], [Djo]) and time scales (see, e.g., [Bhn], [B&A1], [B&P2]) counterparts – extension to wavelets is demonstrated in [DPC-M5]. A common feature of all these NIDs is that they are built on the foundation of the step size as the correct denominator for the measure of rate of change, either directly (local NIDs) or indirectly (non-local NIDs), a paradigm which the ESDDFD refutes. The ESDDFD approach recovers all the NIDs referred to above as limits of the exact representations described hereafter.

The methods discussed in this Thesis do not apply to NID formulations such as the Grünwald–Letnikov (GL) and Kolwanker-Gangal (KG) [K&G1] type (see [Cpls]for a summary list of NIDs), which use the expression $(t - t_0)^\alpha$ or $h^\alpha$ as an acceptable denominator for the derivative. It is suggested in [DPC-M3] and [DPC-M4] that such constructions may have a fatal flaw: they correspond to no exact derivative representation governed by SSP(G) for $\alpha \neq 1$ and therefore cannot be used to model phenomena governed by those rules for $\alpha \neq 1$. Since waves are governed by SSP(G) at spectral level, those NIDs cannot be used to accurately model wave phenomena. This view mistakes linear or power law wave signatures for linear or power law waves; that is, while the Debye and KWW signatures (2.1.1) - (2.1.4) have linear or power law behavior, their sources are not linear or power law waves.

The basic challenge with most differential equations models is that they lack analytical solutions, so that the properties they predict are mostly obtained through mathematical analyses of various asymptotic approximations and numerical simulations. The first step for mathematical analysis is to replace the difference ADRE (3.2) with the differential equation ADRE (3.1) by taking the limit as the time and space steps approach zero, thereby insinuating the time and space steps as integral components of the model as denominators for all derivatives. There are two main problems with this difference-differential replacement strategy. The first problem is that all phenomena observed through the limiting differential equations model have now been labelled as the same and any ability to recover each one as distinct from all others (and remember, there is a whole ocean of them!) is lost. The second problem is the mistaken assumption that 'limit approaching zero' and 'equal to zero' are interchangeable (Fig 2 shows that the two are not interchangeable).

**4 Discrete Modeling DE Models of ADR Brownian Motion and Sub-diffusion Processes**

The motivational problem considered in this thesis is therefore the construction of the discrete analog of the ADRE, subject to various initial and boundary conditions, that faithfully replicates specified properties of the ADRE. The thesis offered here is that such construction must be based on the exact discrete representation of the derivative derived from the SSP at spectral level; any other construction will not have both the simplicity and accuracy of the ESDDFD method for numerically modeling wave behavior down to the spectral level.

The simplest way to construct a discrete model from a continuous DE is using the finite difference method, whose first goal, as is that of any DE discretization method, is the systematic discretization of the derivative, that is, its representation on a discrete lattice. The goal is therefore to construct a strategy for replacing each of $H, W, \frac{\partial W}{\partial t}, \frac{\partial W}{\partial x}, \frac{\partial^2 W}{\partial x^2}$ by its discrete counterpart in such a way that the resulting discrete model retains fidelity to the original model. To enable this, the traditional approach (which we use for our method) is to start (see, e.g., [L&W] Section 12.1) by assuming that the rectangle $\Omega \equiv [a,b]\times[c,d]$, where $(t,x) \in \Omega$, is replaced by the discrete lattice $\{t_n := t_0 + n\Delta t, \ x_m := x_0 + m\Delta x \}_{m,n\geq 0}$, where $(t_0, x_0) = (a,c)$, and $\Delta t, \Delta x > 0$ are, respectively the time and space step sizes. The approximation of the solution of (3.5) at point $(t_n, x_m)$, is then to be obtained as a sequence $\{W_m^n \approx W(t_n, x_m)\}$ that solves the discrete model in the (implicit or explicit) form,

$$W_m^{n+1} = F(t_n, x_m; W_m^n, W_m^{n+1}, W_{m+1}^n, W_{m-1}^n). \tag{4.1}$$

While there are many methods for the numerical approximation of IVP for (3.5), they are NID definition-dependent and none exists that is parallel to the approach proposed here, one enabled by exact difference quotient representations of the derivative from Fourier-Laplace transform space and applicable across a wide range of settings and applications. The aim of this Thesis is therefore to propose a universal representation of the NID that is motivated by the discrete modeling of non-integer-order differential equations (NIDEs), that is based on the fundamental SSP in Fourier-Laplace transform space, and requires no knowledge of the mathematical theory of NIDs.

### 4.1 The Traditional Standard (Euler) Method

Once the lattice is defined, the traditional standard way is to start with the assumption that the standard Euler scheme,

$$H\left(W_m^n, \frac{W_m^{n+1}-W_m^n}{\Delta t}, \frac{W_{m+1}^n-W_m^n}{\Delta x}, \frac{W_{m+1}^n-2W_m^n+W_{m-1}^n}{(\Delta x)^2}\right) = 0, \tag{4.1.1}$$

is a good approximation of (3.4) for very small time and space step sizes, which is justifiable as first and second order truncations of a Taylor series for the first and second derivatives. Improvements are then made by some systematic manipulation of this starting ansatz using the geometry of the problem or some clever mathematical construct designed to preserve some property of the original system.

A fundamental flaw in this approach is that each well-posed boundary/initial value problem for a single ADRE, known to result in a single dynamical system, results in an infinite number of dynamical systems just by varying choices of $\Delta t$ and $\Delta x$ with all other parameters fixed. This appears to be an unfortunate intersection of the rate of change and the slope of the line tangent a given curve; the Taylor series truncation used in deriving ADRE (3.1) leads to equating 'limit approaching zero' with 'being equal to zero'. It is like viewing, in Figure 2 below, the multi-signals in b) through the single lens a) that assumes that they are all the same because they appear the same; however, c) shows that the multi-signals never actually become the same no matter how close to the origin you get.

While the single lens view may yield very good information about the common features of all the wave signatures that look the same, there is no way to distinguish the individual signatures. In a universe where precision is of paramount importance, this inability to distinguish very similar and yet distinct signatures is a fundamental flaw in any numerical method built on such a foundation; the spectral space foundation of the proposed method removes this flaw.

In current NIDE discretization methods, the first step is the expression of the fractional derivative in terms of the corresponding integral representation in space-time (see, e.g., [Diet],[ Chetal],[Ren], [Brck],[A&B],[Dja], [Grrp], [K&B]). The result of the integration representation is that the resulting discrete model has the time and space time-steps insinuated as the foundation derivative descriptors, thereby losing all spectral information from all derivative representations.

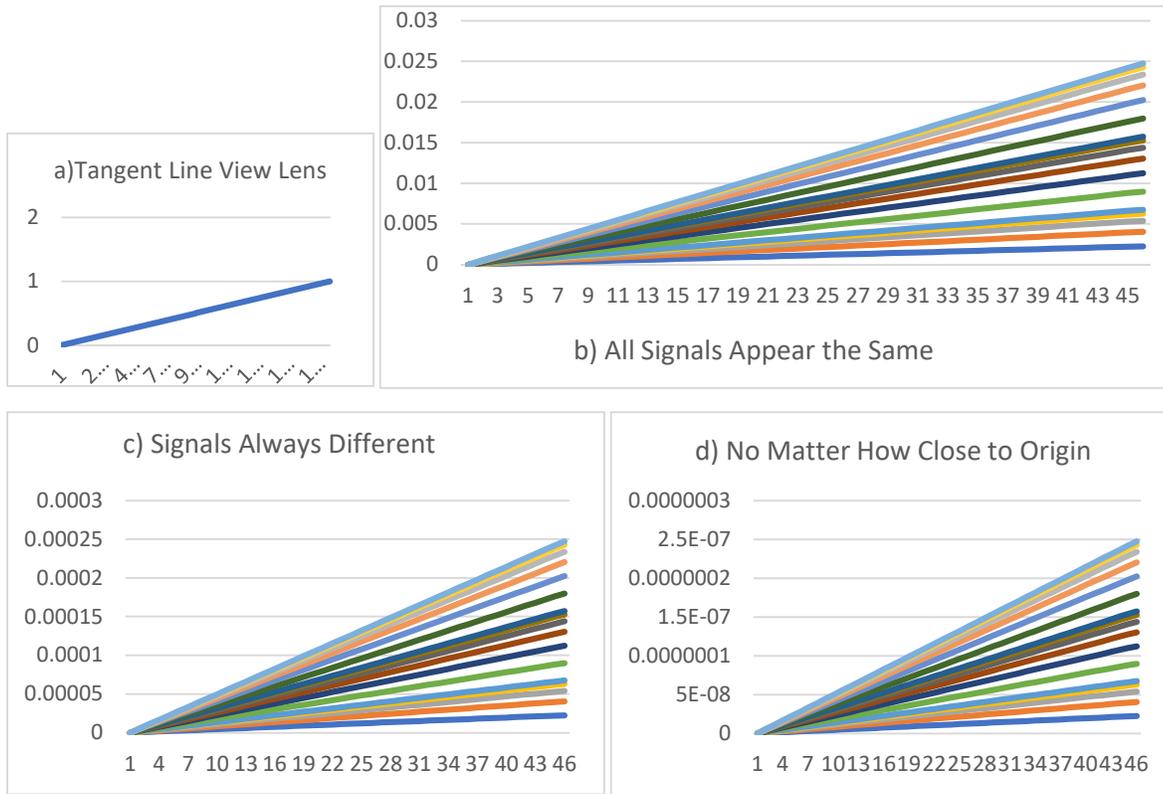

Figure 2 (Traditional Euler Lens Shop)    Distorted views; with the traditional denominator (h) the lines in (b) look the same as the one line in (a) when step size is below h=0.3 and are subsequently viewed as being the same. However, (c) and (d) show that no matter how close to zero one gets, different lines remain different.

In our view lens paradigm, current NIDE discretization methods design lenses that are excellent at approximating most signatures, but they cannot distinguish them from each other at spectral level. It is akin to having very flexible lenses for the 'line view' of Figure 2.

## 4.2 The Non-standard Finite Difference (NSFD) Method

The (Mickens) NSFD method was designed to specifically address the problem that the Euler method defines different dynamical systems for the same IV/BVP for different step sizes. Its solution to this problem is the NSFD denominator and non-linear term discretization with positivity construction safeguards.

In the late 1980's Mickens realized that there was nothing particularly special about the traditional step size denominator; any function with the same limit could serve the same purpose. This was the beginning of the traditional nonstandard finite difference (NSFD) method, which is summarized in [Mick1]. The NSFD way is to start with the assumption that the standard Euler scheme,

$$H\left(W_m^n, \frac{W_m^{n+1}-W_m^n}{\phi_1(\Delta t)}, \frac{W_{m+1}^n-W_m^n}{\phi_2(\Delta x)}, \frac{W_{m+1}^n-2W_m^n+W_{m-1}^n}{(\psi(\Delta x))^2}\right) = 0, \tag{4.2.1}$$

where each of $\phi_1, \phi_2, \psi$ have the same limit behavior, as the step size is a good approximation. Moreover, these denominators can be obtained in various prescribed ways (see, e.g. [Mick1], [Mick2], [A&L],[JS], [Wds], [WdsTh])) such that, with prescribed treatment of the non-derivative terms, the resulting discrete models have the same behavior as the model for all step sizes, thereby describing the same dynamical system for all step sizes. The denominator connection to stability was made mathematically more precise with the notion of elementary stability of discrete models (i.e. having the same behavior as the continuous model for all step sizes) in the work of Lubuma and his collaborators, which is summarized in [A&L]. That theory establishes as elementary stable all models with denominators constructed according to rules that include application of SSP(L) in ordinary space-time.

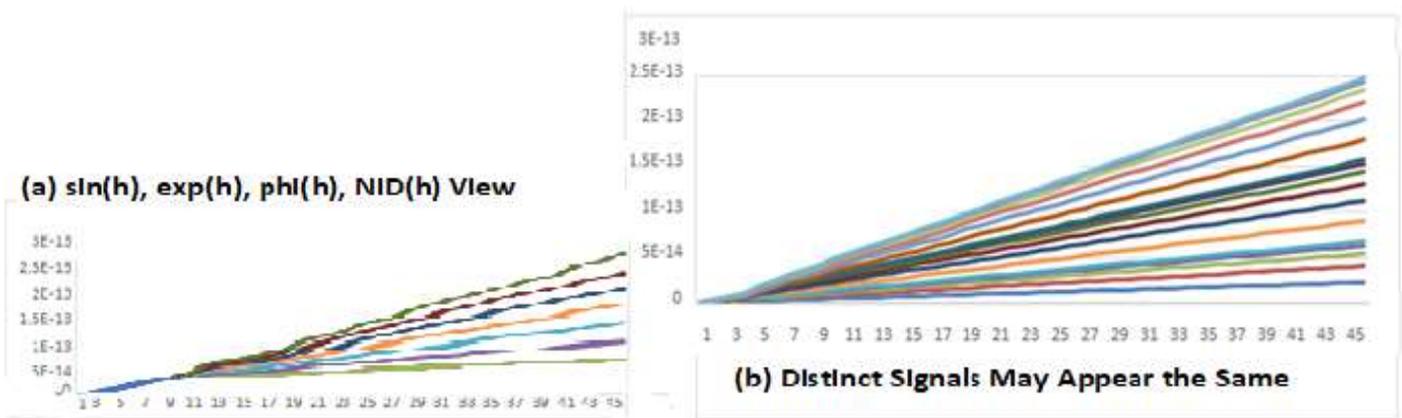

Figure 3 (Mickens NSFD-NID-Euler Lens Shop)    Better than Traditional, but still distorted views; even with the NSFD denominator $\phi(h)$ in NIDEs, the lines in (b) <u>look the same</u> as the ones line in (a). However, besides there being a small number of distinct views, those views include lines whose signature is not that of waves with any known behavior.

While it partially solves the problem of discretization generating different dynamical systems from the same model by mere change of step size, the Mickens NSFD methodology still has a foundational flaw: spectral data exclusion from the various prescribed ways of denominator construction. The method starts with the assumption that the limit denominator, the step function, is acceptable, thereby classifying all same-looking signatures with the same limit behavior as the same without the ability to correctly ascribe any described behavior to the correct source out of those same-looking signatures. Acceptable denominators in the Mickens method therefore include the following functions:

$$\phi(h) = h, 1 - e^{-h}, e^h - 1, \sin(h), 1 - \frac{h+1}{h}, \tag{4.2.2}$$

all of which are fundamentally flawed; further discussion about these denominators, the order of numerical schemes, and the answer to Mickens' original question are in Chapter 2 of the Thesis. The proposed method overcomes this flaw by requiring that the Fourier-Laplace (or their counterparts in different systems) transform space be the native discretization space for all ADREs and their generalizations; that is, the fundamental SSP must be applied in spectral transform space.

A further constraint for the Mickens NSFD method is using as one of its core tenets is the discretization of the SSP(LI) propagator from its exact solution from a group theoretic property of the exact solution. This enables use of an exact solution operator, $E(t_0, y_0; .)$, of an IVP to obtain the exact discrete model of that IVP using a simple transformation (see, e.g, [Mick1], [A&L], [JS]),

$$y(t) = E(t_0, y_0; t) \rightarrow y_{n+1} = E(t_n, y_n; t_{n+1}), \tag{4.2.3}$$

which is then used to obtain the discrete form of the derivative. That group property is, to my understanding, not available for non-integer solutions.

It appears that recent attempts to extend the NSFD method to fractional DEs have been hampered by an effort to retain this group theoretical property as well as the adherence to the 'step size as basic denominator' (see, e.g., [GnzP], [D&D], [App], [Swe], [Dáv]) as correct. This h-dependence is even reflected in the definition of Mittag-Leffler stability ([Lietal], (see also [Sene], [L&M]). As a result, all these efforts end up with a derivative denominator of the form $(t - t_0)^\alpha$, which is then discretized by Eqn (4.2.3) as $(h)^\alpha$; this is equivalent to truncating the NSFD denominator, a fatal flaw discussed further in Chapter 2 of this thesis.

In our view lens paradigm, current NIDE NSFD discretization methods design lenses that are excellent at approximating most signatures; however, they cannot distinguish them each from the other and their view include signatures that do not represent any waves from spectral space. It is akin to having very flexible lenses with multi-view attachments to the 'line view' lens, as in Figure 3.

## 4.3 The Exact Spectral Derivative Discretization Finite Difference (ESDDFD) Method: Exact Measures for Viewing Wave Signatures

The (Clemence-Mkhope) ESDDFD method is designed to specifically address the twin problems that the Mickens method carries no spectral information about the waves being modeled and that its extension to NIDEs is hampered by group-theoretical considerations, reflected in its h-dependence. The ESDDFD solution to both these problems is the novel derivative denominator construction in Fourier-Laplace transform space, with non-linear term discretization and positivity construction safeguards.

The impetus for the new approach was the observation that behavior of the form classified as Debye and KWW patterns (2.1.1)- (2.1.4) in Section (2.1) is captured by the following local and non-local propagators for SSP(1G), $0 < \alpha \leq 1$, which do not include (2.1.1)- (2.1.4):

Local $\qquad\qquad y(t) = y(t_0)\exp(-\lambda t^\alpha)$ (4.3.1)

General Local $\qquad y(t) = N(t_0,)\exp(\psi(t, \lambda, \alpha))$ (4.3.2)

Non-local $\qquad\quad y(t) = y(t_0)E_\alpha(-\lambda t^\alpha),$ (4.3.3)

General Non-local $\quad y(t) = y(t_0)E_\alpha(\psi(t, \lambda, \alpha))$ (4.3.4)

where is the Euler (exponential) function,

$$\exp(z) = e^z == \sum_{k=0}^{\infty} \frac{z^k}{k!},$$

and $E_\alpha(z)$ is the one-parameter Mittag-Leffler (ML) function,

$$E_\alpha(z) = E_{\alpha,1}(z) = \sum_{k=0}^{\infty} \frac{z^k}{\Gamma(\alpha k+1)} \approx 1 + \frac{(z)}{\Gamma(\alpha+1)} + \cdots = 1 + \frac{z}{\Gamma(\alpha+1)},$$

where $E_{\alpha,\beta}(z)$ is the two-parameter ML function,

$$E_{\alpha,\beta}(z) = \sum_{k=0}^{\infty} \frac{z^k}{\Gamma(\alpha k+\beta)}.$$

The functions (4.3..1)-(4.3.4) serve as the basic building blocks for the proposed non-limit derivative with a universally defined denominator measure specified in terms of these blocks; they are the probability distribution functions (pdf) of 'pure', spectral waves (see [Fyn] and [B&CC] for respective physicist and biomathematical pdf perspectives). While the asymptotic behavior of the functions (4.3.1)-(4.3.4) is linear or power law, none of them are the linear or power law functions (2.1.1)-(2.1.4), which do not behave according to the SSP and therefore do not represent waves. The approach proposed here is to replace any advection-diffusion-reaction differential equation by its difference quotient equivalent, but with the standard time and space step size or Mickens NSFD denominator quotients replaced with ones with an ESDDFD denominator specifically obtained from Fourier-Laplace transform spaces to preserve all spectral and other behavior as follows:

- The time evolution measure is obtained from the exact discretization of the Fourier transform of the reactionary heat sub-equation.
- The diffusion measure is obtained from the exact discretization of the Laplace transform of the transient diffusion-reaction sub-equation.
- The advection measure is obtained from the exact discretization of the (once-integrated) steady state advection-diffusion sub-equation.
- All non-linear terms are discretized non-locally according to the Mickens NSFD rules.
- All terms are discretized to preserve positivity according to the Mickens NSFD rules.

The resulting numerical models are all unconditionally convergent and carry complete (as provided by the model):

- spectral (for all wave modes), diffusion, and reaction information in the time evolution denominator, obtained from SSP(I)
- spectral (Laplace), diffusion, and reaction rate information in the diffusion denominator, obtained from SSP(II)
- advection and diffusion information in the advection denominator, obtained from SSP(I)
- dynamical system properties of the continuous model – and different time or space step sizes do not define different dynamical systems.

While the exposition given here is for deterministic models and focuses on one space dimension, the following is provided in [DPC-M4] as evidence for extensions as follows:

- a template for extension to two-dimensional and three-dimensional Laplacians is given,
- the Black-Scholes [Demp] equation is given as evidence of possible extension to stochastic processes,
- SSP(I) and SSP(II) measures for local Debye processes on time scales ([B&P]) are given,
- extension paths for wavelet and string theories are suggested.

The discrete models of differential equations resulting from the approach proposed here will enable simple finite difference simulation of a wide variety of phenomena that offers a more refined, and more reliable, views of phenomena across many spheres of human endeavor than is currently possible from any methods built on any limiting definition of the derivative. In our world of wave view windows, the view from the resulting simulations may be about as close as we get to each ray having its own view lens from any window, so that 'infinitesimally similar' objects will look different if they indeed are, such as the 'wave lines' in Figure 4 below.

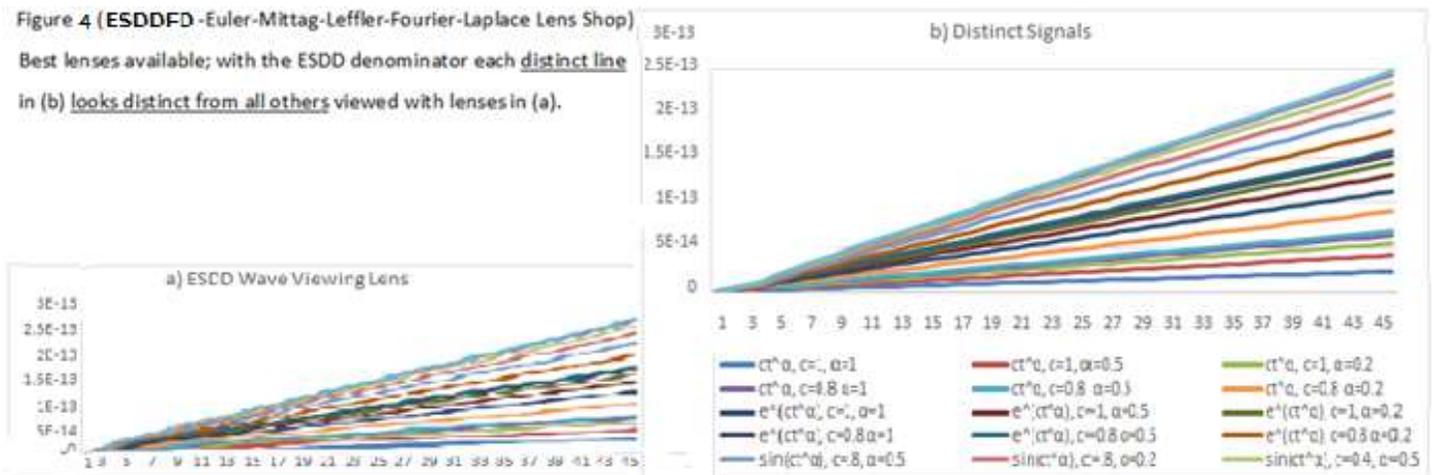

Figure 4 (ESDDFD -Euler-Mittag-Leffler-Fourier-Laplace Lens Shop) Best lenses available; with the ESDD denominator each distinct line in (b) looks distinct from all others viewed with lenses in (a).

The ESDDFD method diverges from current NIDE discretization methods by replacing the first step in those methods by transformation into Fourier-Laplace spectral space. The result of the spectral transformation is that the resulting ESDDFD models have denominators that incorporates all spectral information into their derivative representations. The new ESDDFD-Mittag-Leffler-Fourier-Laplace Lens Shop template produces a unique lens for each signature by specification of wave spectral data. Each ray in Fig 3 has its own view lens.

## 5 Example Contrasting Traditional Euler, Mickens NSFD, and Spectral ESDDFD Methods

To compare the methods discussed, Eqn (5.1) below, (Eqn (7.3.4) in [Mick1]) is discretized using each of the standard Euler, NSFD, and ESDDFD methods.

$$\frac{\partial u}{\partial t} = a \frac{\partial^2 u}{\partial x^2} + bu \tag{5.1}$$

where, $a, b$ are constants with $a \geq 0$. It is to be noted modeling requirements that must be satisfied for all methods are not a step size restriction (see, e.g., [Mick1], [M&K]).

### 5.1 The Traditional Standard (Euler) Method

The traditional standard (Euler) method for (5.1), given in Eqn. (5.1.1) below, is based on the time and step sizes being the derivative denominators, which is justified by the defining the derivative as the slope of a line tangent to a curve and then treating DE models as truncations of Taylor expansions.

$$\frac{u_m^{n+1} - u_m^n}{\Delta t} = a \frac{u_{m+1}^n - 2u_m^n - u_{m-1}^n}{(\Delta x)^2} + bu_m^n \tag{5.1.1}$$

While accurate (and can be made more accurate with the improved Euler method) for very small step sizes (and so can serve as a benchmark method), it performs poorly in comparison to other standard (i.e., based on step size denominators) methods to be useful in practice (the smaller the step sizes, the more expensive the method due to implementation time) and is poor at simulating various types of initial and boundary values encountered in applications. Further, the method defines different dynamical systems for the same IV/BVP for different step sizes, with some exhibiting expected dynamical behavior and some producing chaotic behavior for non-chaotic phenomena. Moreover, the method carries neither diffusion and reaction effects in its denominators nor spectral information in its construction; therefore it, or any method with it as a foundation, cannot accurately simulate any diffusion and reaction effects or spectral behavior of the phenomena it models.

### 5.2 The (Mickens) NSFD Method

The NSFD model for (5.1), given in [Mick1] and as (5.2) below, is based on the exact discretization as follows of sub-equations of (5.1), which are ODEs that embody SSP(L) in ordinary, that is, $(x, t)$, space:

$$\frac{du}{dt} = bu \tag{5.2.1}$$

$$a\frac{d^2u}{dx^2} + bu = 0 \tag{5.2.2}$$

The exact discrete representations of (5.2.1) and (5.2.2) are, respectively,

$$\frac{u^{n+1}-u^n}{\phi(\Delta t,b)} = bu^n \tag{5.2.3}$$

$$a\left(\frac{u_{m+1}-2u_m-u_{m-1}}{(\psi(\Delta x,b/a))^2}\right) + bu_m = 0, \tag{5.2.4}$$

where the NSFD denominators are given by

$$\phi(\Delta t,b) = \frac{1}{b}(e^{b\Delta t}-1) \text{ and } \psi(\Delta x,b/a) = (2\sqrt{a/b})\sin\left(\sqrt{b/a}\,\frac{\Delta x}{2}\right). \tag{5.2.5}$$

An explicit scheme (an implicit one is similarly obtained) for (5.1) is then obtained, by combining (5.2.3) and (5.2.4), as (5.2) below:

$$\frac{u_m^{n+1}-u_m^n}{\phi(\Delta t,b)} = a\left(\frac{u_m^{n+1}-2u_m^n-u_{m-1}^n}{(\psi(\Delta x,b/a))^2}\right) + bu_m^n, \tag{5.2}$$

While it carries the diffusion and reaction effects in its denominators and is competitive with or better than most models from existing discrete methods in its simplicity and implementation, the model (5.2) carries no spectral information about the waves being modeled; therefore it, or any method with it as a foundation, cannot accurately simulate any spectral behavior of the phenomena it models. Moreover, the method's fundamental dependence of exact solutions in time-space appears to create a hurdle in its extension to NIDEs due to group-theoretical considerations.

### 5.3 The (Clemence-Mkhope) ESDDFD Method

The ESDDFD method for (5.1) is based on the exact discretization as follows in Fourier and Laplace transform spaces of the complete Eqn (5.1), which does not embody (SSP(L) in ordinary $(x,t)$ space. However, in Fourier transform space, with $U(k,t) = \mathcal{F}_x\{u(x,t)\}$, where $k$ is the Fourier spectral parameter (wave modes), (5.1) becomes SSP(LI):

$$\frac{dU}{dt} = (b-ak^2)U \tag{5.3.1}$$

whose exact discrete representation is SSP(GI),

$$\frac{U^{n+1}-U^n}{\phi(\Delta t,b,k)} = (b-ak^2)U^n, \tag{5.3.2}$$

where the ESDDFD denominator is given by

$$\phi(\Delta t,b,k) = \frac{1}{(b-ak^2)}\left(e^{(b-ak^2)\Delta t}-1\right). \tag{5.3.3}$$

Similarly, in Laplace transform space, with $Y(x,s) = \mathcal{L}_t\{u(x,t)\}$, where $s$ is the Laplace spectral parameter (Laplace wave modes), (5.1) has the same linear form as SSP(LII):

$$\frac{d^2Y}{dx^2} + ((b-s)/a)Y + \frac{1}{a}u(x,0) = 0 \tag{5.3.4}$$

whose exact discrete representation is SSP(GII),

$$\left(\frac{Y_{m+1} - 2Y_m - Y_{m-1}}{(\psi(\Delta x, b/a, s))^2}\right) + ((b-s)/a)Y_m + \frac{1}{a}u(x_m, 0) = 0, \tag{5.3.5}$$

where the ESDD denominator is given by

$$\psi(\Delta x, b/a, s) = \left(2\sqrt{a/(b-s)}\right)\sin\left(\sqrt{(b-s)/a}\,\frac{\Delta x}{2}\right). \tag{5.2.5}$$

An explicit scheme (an implicit one is similarly obtained) for (5.1) is then obtained, by using (5.3.3) and (5.3.6) as exact derivative representations, as (5.2) below:

$$\frac{u_m^{n+1} - u_m^n}{\phi(\Delta t, b, k)} = a\left(\frac{u_m^{n+1} - 2u_m^n - u_{m-1}^n}{(\psi(\Delta x, b/a, s))^2}\right) + bu_m^n, \tag{5.3}$$

The model (5.3) carries the diffusion and reaction effects as well as spectral information about the waves being modeled in its denominators; Fourier spectral information is in the time evolution denominator while Laplace spectral information is in the space diffusion denominator. Due to the diffusionless case of (5.1) having the same SSP (I) form in ordinary and Fourier transform spaces and its steady state case having the same SSP (II) form in ordinary and Laplace transform spaces, the NSFD model is recovered as a special case of the ESDDFD model in these special cases.

Based on its construction, it is our claim that the ESDDFD method is competitive with or better than most (if not all) models from existing discrete methods, in its simplicity and implementation, and will outperform those methods in simulating any spectral behavior of the waves being modeled.

*6 Twelve Mathematical/Computational Consequences of Proposed Universal Wave View Modeling*

The SSP-based, ESDDFD approach has several, significant, short term and long-term mathematical and computational consequences; some, outlined below, are the topics of the remaining Thesis chapters.

1) Chapter 2 – Order and Truncation. In answer to the original motivating Mickens question, the definition of order of numerical method is revisited (and revised) – truncation has a new meaning. Key points in this chapter are as follows: (a) The meaning of the order of a numerical scheme is re-defined, with its order referring to spectral parameter order instead of that of the step size. (b) The traditional step size denominator is catastrophically unsuited to 'truncated exponential denominators' (e.g., at unit step size, odd truncation results in non-existent derivatives, while even truncation does not).

2) Chapter 2 – Elementary Stability. While some bounds still need to be calculated as per model and/or NID (e.g., fractional models converge in wider angled regions that do integer ones, but the angle depends on the NID), some arguments about elementary stability (which are about specifying classes of similarly behaved dynamical systems) may become moot since our schemes converge unconditionally according to the theory of [A&L] for the NSFD method interpreted through the correct viewer, except for step size relationships imposed by physical modeling constraints.

3) Chapter 3 [DPC-M3]– Core Theory. All non-integer derivatives with solutions to the Cauchy IVP for SSP(LI) are limit cases for their SSP(GI) counterparts. Key points in this chapter section are: (a) from the way NIDEs arise in modeling, only Caputo type derivatives fall in this category; (b) the Riemann-Liouville and Caputo versions of the NIDEs describe the same phenomena from different viewpoints – the relationship between the two is 're-discovered' from modeling considerations and the RL memory effect is modeled without boundary or initial value issues, suggesting that the issue with such problems may be due to looking at un-justifiable conditions or models. For example, (i) (see [M&K], Section 5.1), while non-trivial IV/BVPs may be posed for ${}_0^C D_t^\alpha u(x,t) = a\frac{\partial^2 u}{\partial x^2} + bu$, they cannot be posed for ${}_0^{RL}D_t^\alpha = a\frac{\partial^2 u}{\partial x^2} + bu$ – the later models waves $u(x,t)$ whose initial value, $u(x,0)$, is zero while the former models non-trivial waves; (ii) interpreting the condition $u(x,0) \sim Cx$ as meaning that $u(x,0) = Cx$ (as appears to be in [CL&V]) is saying that all functions that look like $Cx$ are the same as $Cx$, which they are not.

4) Chapter 3 [DPC-M 3] – The Local to non-local Switch. (a) The conformable fractional derivative (CFD) is the link between local and non-local NIDs – a result due to Mainardi [FM] to this effect is re-discovered and generalized; methods termed the 'conformable Euler method' [Xin], [Mhmmdz] and 'generalized Euler method' [O&M] are shown to be incorrect, respectively in [DPC-M4] (see also [DPC-M1]) and [DPC-M5], using the CFD. (b) There is a relationship, at least for the limit case, between the local and non-local Cauchy IVP, universally specified in terms of SSP(GI), which leads to a relationship, at least in a limit sense, between the (non-local) M-L function in terms of the (local) Euler function.

5) Chapter 3 [DPCM-3] – Extension Example. An accurate discrete model may be constructed that fully (for all wave modes) reproduces Fourier and Laplace spectral effects, through a consistent application of two self-sameness principles (SSP) in the appropriate spectral transform spaces. The following examples in one space dimension are discretized in [DPCM-3] as a demonstration of the SFDDFD method: (a) the Fokker-Planck equation (FPE) from [M&K] for both Brownian motion – a variety of known FPE properties, such as the Bose-Einstein and Boltzman identities are easily recovered from the ESDD model; (b) the fractional Fokker-Planck equation (FFPE) from [M&K] (see also [Z&K], [K&G2], [K&B], the properties extending those in (a) are easily recovered from the ESDDFD model.

6) Chapter 4 [DPC-M4] - Dimensional and Other Extensions. An accurate discrete model may be constructed (whenever the heat, wave and Laplace equations or their generalizations hold true) that fully (for all wave modes) reproduces Fourier and Laplace spectral effects using the SSP paradigm; the

ESDDFD construction has the same formulation (i.e., is invariant) across NID definitions, dimensions, time scales, deterministic or stochastic (respectively reflecting 'natural' and human-driven processes or random walks), and wave views (i.e., classical, wavelet, or string). The SSP application extension to two and three space dimensions, as well as to some stochastic processes, is presented in [DPC-M6] and extension to time scales, general stochastic processes, and wavelet views are conjectured – constructions are presented with conjectured functions. In 2-D and 3-D space, Debye and KWW time patterns ($\sim Ct^\alpha, 0 < \alpha \leq 1$) as well as space patterns ($\sim Cx^\beta, 0 < \beta \leq 1$) originate in Fourier-Laplace transform space, with both spectral effects (including cross-coordinate axis Fourier effects) being present in both denominators.

7) Chapter 4 [DPC-M6] – More Examples. The accurate discrete modeling of fractional (in space and/or time) advection-diffusion-reaction systems using the SSP paradigm does not require an a-priori notion for the mathematical theory of any derivative; one can construct accurate discrete models of fractional and other derivatives BVPs without knowledge about those derivatives other than that of their existence. The following examples are discretized for $0 < \alpha \leq 1$ as further ESDDFD demonstration: (a) 1-3 dimensional elliptic ADRE with Neumann Conditions as considered in [Nkt]; (b) the parabolic ADRE with Wentzell-Robin Conditions as considered in [Nkt], (c) the Black-Scholes equation for the variety of options described in [Demp]; (d) the time scales ([B&P1]) ODE considered in [Bhn]. The existence, uniqueness, regularity, and continuous dependence on the data of the resulting ESDDFD models will be discussed in future articles.

Chapter 5 – Conjectures and Other Consequences. The last chapter of the Thesis collects various conjectures and outlines propositions for the ideas listed as (8) – (12) below.

8) It is conjectured that a conundrum for infinitesimal Calculus (see, e.g., [K&T] or [WT for a quick introduction) may be resolved: that of assuming a non-zero denominator which is then arbitrarily set equal to zero at evaluation [WT] – the issue does not arise if a probability measure is used as the infinitesimal derivative denominator.

9) An interesting curiosity follows for group theory: The Mickens NSFD method was not extended to NIDEs mainly due to the belief that the exact solution to exact scheme was dependent on some dynamical sub-group property. However, the construction for the exact representation works the same for integer and NID SSP laws. Therefore, perhaps re-viewing the dynamical systems for NIDE IV/BVPs through ESDDFD lenses may reveal something new about group structures at spectral level?

10) An interesting possibility follows for number theory: Through logarithmic scaling – an immediate consequence of the exponential function – and the following is conjectured: The number of primes in the interval $(1, \ exp(\pi^{2n}))$ can be counted through a logarithmic map of intervals; that is, we conjecture that there is a Collatz type logarithmic spiral (see., e.g. [Reid21]) mapping from the number of primes in $(1, \ exp(\pi^{2n}))$ to those in $(0, \ \pi^{2n})$. A related question is how many primes there are in $(\pi^{2n}, \ 10^n)$, where $n$ is a positive integer, a curious relationship between $\pi$ and 10. The Collatz spiral argument (which, in my opinion, proves the Collatz conjecture (see, e.g., [Tao], [Reid20a])) is presented in [Reid20b] (see also [Reid20a]); part of that argument is now available as [Reid21]. An observation of Reid's argument is that in the spiral function, $r = 2 * 2^{\theta/\pi}$, $\pi$ appears to be no coincidence (nor is the fact that $\pi^3 \approx 3\pi^2 \approx e\pi^2$): ours appears to be an "$e - \pi$" universe'! Moreover, since the Collatz conjecture is an ordered pair, $(2, 3)$, construction (as is Reid's), are there other 'Collatz pairs'? If there are, that would imply that no matter where you begin, there is a 'binary prime path', that leads 'home' – how long it takes just depends on your prime pair.

11) An interesting possibility follows for number theory: By considering Fermat triplets and using a logarithmic map argument, Fermat's last theorem (see, e.g., [MD]) may be a simple special case of the statement that you can only inscribe up to (n+1)-sided integer polygons inside an n-dimensional sphere; that is, you cannot have an (n+2)- gon inside an n-sphere whose sides are integers – the only integer-sided polygons are those with their number of sides up to one larger that the sphere dimensions. Fermat's theorem is that only up to a triangle can be inscribed inside a circle, and not a rectangle or any polygon with number of sides larger than three; if that is so, then perhaps Beal's conjecture [MD] is not a 'correct' generalization of Fermat's last theorem?

12) Immediate computational advantages:

a) One needs no a-priori knowledge about the mathematical theory of fractional or fractal derivatives, in any time scales, dimensions, or wave views, to construct and execute models in those paradigms.

b) Any conditionally stable algorithm based on the traditional limit definition of the derivative can be rendered unconditionally stable as well as more accurate with the replacement of the traditional Euler or NSFD denominator by one derived from the ESDDFD method. While unconditional stability can be achieved from the NSFD denominator, that method lacks the systemic incorporation of spectral measures into derivative denominator construction.

**Funding:** This research did not receive any grant from funding agencies in the public, commercial, or not-for-profit sectors.

**Acknowledgements**
I wish to thank, above all, my friend and life partner, Belinda G. Brewster-Clemence (Mkhope), who listened to my musings (some of these were first thing in the morning – before she even had a chance for coffee) over the years about everything contained herein, read countless drafts of this manuscript, and is responsible for all the graphs included herein. (Further acknowledgements are in the Thesis.)

**References**

[Mick 1] Mickens, R.E. (1994). *Nonstandard finite difference models of differential equations*.

[Fyn] Feynman, R.P., Leighton, R.B., and Sands, M. (1977, Sixth Printing). *The Feynman Lectures on Physics*, Volume I. Addison-Wesley Publishing Co

[HS] Hinton, D., Shaw, J.K (1981). On Titchmarsh-Weyl M(λ)-functions for linear Hamiltonian system. *Journal of Differential Equations,* 40: 316–342
DO: 10.1016/0022-0396(81)90002-4

[HKS], Hinton, D., Klaus, M., and Shaw, J. (1988). Levinson's theorem and Titchmarsh-Weyl m(λ) theory for Dirac systems. *Proceedings of the Royal Society of Edinburgh: Section A* Mathematics, (109). DO: 10.1017/S0308210500026743

[DPC94] Clemence, D. P. (1994). Titchmarsh-Weyl theory and Levinson's theorem for Dirac operators. *Journal of Physics A*, 27:7835-7842

[C&L] Coddington, Earl A. and Levinson, Norman (1984 UK ed. Edition). *Theory of Ordinary Differential Equations*. Malabar, Fla.: R.E. Krieger, 1984. Originally published: New York: McGraw-Hill, 1955. (International series in Pure and Applied Mathematics)

[Khl], Khalil, R., Al Horani, M., Yousef, A. and Sababheh, M. (2014). A new definition of fractional derivative. *Journal of Computational and Applied Mathematics*, 264: 65-70.

[Musr] Musraini M., Rustam Efendi, Endang Lily, Ponco Hidayah (2019). Classical Properties on Conformable Fractional Calculus. *Pure and Applied Mathematics Journal*, 8(5): 83-87

[Vld] Valdes, Juan E. Nápoles, Guzmán Paulo M., Lugo, Luciano M., and Kashuri, Artion (2020). The local Generalized Local Derivative and Mittag-Leffler Function. *Sigma J Eng & Nat Sci*, 38 (2): 1007-1017.

[KCi] Karcı, A. (2013). A New Approach for Fractional Order Derivative and Its Applications. *Universal Journal of Engineering Science*, 1: 110-117.

[Z&LC] D Zhao, M Luo-Calcolo (2017). General conformable fractional derivative and its physical interpretation. 2017 - Springer

[Almd] Almeida, Ricardo, Małgorzata Guzowska, and Tatiana Odzijewicz (2016). A remark on local fractional calculus and ordinary derivatives. *Open Math.* 2016; 14: 1122–1124;

[Imb] Imbert, Alberto F. (2019). *Contributions to Conformable and Non-Conformable Calculus*

[GJ] G. Jumarie (2006). New stochastic fractional models of the Malthusian growth, the Poissonian birth process and optimal management of populations. *Math. Comput. Modell*. 44: 231-254.

[Cpls] Edmundo Capelas de Oliveira1 and José António Tenreiro Machado (2014). A Review of Definitions for Fractional Derivatives and Integral. *Mathematical Problems in Engineering* Volume 2014, Article ID 238459; http://dx.doi.org/10.1155/2014/238459

[TEV1]) Tarasov, Vasily E.  (2016). Local Fractional Derivatives of Differentiable Functions are Integer-order Derivatives or Zero. *International Journal of Applied and Computational Mathematics*, 2 (2): 195-201. DOI: 10.1007/s40819-015-0054-6

[TEV2] Tarasov, Vasily E.  (2018). No Nonlocality. No Fractional Derivative; *Communications in Nonlinear Science and Numerical Simulation*, 62: 157-163. DOI: 10.1016/j.cnsns.2018.02.019

[DPC-M1] (2021 Preprint).  Conformable Euler's Finite Difference Method - Comment https://arxiv.org/abs/2105.10385

[DPC-M 4] Clemence-Mkhope, D.P. (2021 Preprint): On the Conformable Exponential Function and the Conformable Euler's Finite Difference Method

[M&K] Metzler, Ralf and Joseph Klafter (2000). The random walk's guide to anomalous diffusion: a fractional dynamics approach. *Physics Reports* 339 (2000) (1)77. DOI: 10.1016/s0370-15730000070-3. Downloaded 04/29/2021 from:  https://www.tau.ac.il/~klafter1/258.pdf

[B&B] Baldock, G.R. and Bridgeman, T. (1981). *Mathematical theory of wave motion*. Available at: https://ui.adsabs.harvard.edu/abs/1981ehlh.bookQ....B {Provided by the SAO/NASA Astrophysics Data System}

[TEV3] [Tarasov, V. E. (2011). *Fractional Dynamics: Applications of Fractional Calculus to Dynamics of Particles, Fields and Media*, Nonlinear Physical Science, Springer, Heidelberg, Germany.

[Chetal] Chen, W., Zhang, X. D., and Korosak, D. (2010). Investigation on Fractional and Fractal Derivative Relaxation Oscillation Models. *Int. J. Nonlin. Sci. Num*., 11 (2): 3-9

[Z&C] Zill and Cullen. Differential Equations with Boundary Value Problems

[W&W] Weiss G., Wilson E.N. (2001). The Mathematical Theory of Wavelets. In: Byrnes J.S. (eds) *Twentieth Century Harmonic Analysis — A Celebration*. NATO Science Series (Series II: Mathematics, Physics and Chemistry), vol 33. Springer, Dordrecht. https://doi.org/10.1007/978-94-010-0662-0_15

[Grnf] Gorenflo, Rudolf, Francesco Mainardi, and Sergei Rogosin (2019). Mittag-Leffler function: properties and applications. In *Handbook of Fractional Calculus with Applications*, Vol. 1 Basic Theory (A. Kochubei, Yu.Luchko (Eds.) pp. 269-296. (8 Volume Series Ed., J. A. Tenreiro Machado).

[Gmz] José Francisco Gómez, Lizeth Torres, Ricardo Fabricio Escobar (2019). Fractional Derivatives with Mittag-Leffler Kernel. *Trends and Applications in Science and Engineering*, Springer, Feb 13, 2019.

[Pdlb] Igor Podlubny (1999). *Fractional Differential Equations*. Academic Press.

[Ch] Chen, W. (2006). Time-space fabric underlying anomalous diff.usion *Chaos, Solitons and Fractals*, 28 (4): 923-929.

[He] He, J.H. (2011). A New Fractal Derivation. *Thermal Science*, 15 (Suppl. 1): S145-S147

[DPC-M 3] Clemence-Mkhope, D. P. (2021, Preprint). The Exact Spectral Derivative Discretization Finite Difference (ESDDFD) Method for Wave Models: Theory and Examples

[JDL] J. D. Lambert (1991). *Numerical methods for ordinary differential systems*, Wiley, New York.

[L&W] Linz, P. and Wang, R. (2003). *Exploring Numerical Methods: An Introduction to Scientific Computing Using MATLAB*. Jones and Bartlett Learning

[C-F] Caputo, M., Fabrizio M., A New Definition of Fractional Derivative Without Singular Kernel. *Progress in Fractional Differentiation and Applications*, 1(2015),2, pp. 73-85

[A&B] Atangana, A., Baleanu, D. (2016). New Fractional Derivatives with Nonlocal and Non-Singular Kernel: Theory and Application to Heat Transfer Model. *Therm Sci*. 20(2), pp. 763-769.

[Diet] K. Diethelm, N. J. Ford, A. D. Freed, and Yu. Luchko (2005). Algorithms for the fractional calculus: a selection of numerical methods. *Comput. Methods Appl. Mech. Engrg*. 194:743–773.

[Ren] Ren, J., Zhi-Zhong Sun, and X. Zhao (2013). Compact difference scheme for the fractional sub-diffusion equation with Neumann boundary conditions. *J. Comput. Phys.,* 232: 456-467.

[Brck] Brociek, R. (2014). Implicite finite difference method for time fractional heat equation with mixed boundary conditions. *Zeszyty Naukowe Politechniki Śląskiej, Matematyka Stosowana*, 4: 73-87

[Dja] Djida, J.D., I. Area, and A. Atangana (2017). Numerical Computation of a Fractional Derivative with Non-local and Non-Singular Kernel. arXiv: 1610.0717v1 [math.AP] 23 Oct 2016

[Grrp] Garrappa, Roberto (2018). Numerical Solution of Fractional Differential Equations: A Survey and Software Tutorial. Reprinted from: *Mathematics* 2018, 6, 16, doi: 10.3390/math6020016

[K&B] Kumar, S. and Baleanu, D (2020). A New Numerical Method for Time Fractional Non-Linear Sharma-Tasso-Oliver Equation and Klein-Gordon Equation with Exponential Kernel Law. *Front. Phys*. 8:136. Doi: 10.3389/fphy.2020.00136

[Hnsl] E. Heinsalu, M. Patriarca, I. Goychuk, G. Schmid, and P. H¨anggi (2007). Fractional Fokker Plank dynamics: Numerical algorithm and simulation. *Phys. Rev.* E 73, 046133.

[Mick 2]= R.E. Mickens (2007). Calculation of denominator functions for NSFD scheme for differential equations satisfying a positivity condition. *Numerical Methods for PDEs,* Vol. 23.

[A&L] R. Anguelov and J. M-S Lubuma (2001). Contributions to the mathematics of the nonstandard finite difference method and applications. *Numerical Methods for PDE*, 17: 518-543.

[JS] Sunday, J. (2010). On Exact Finite-Difference Scheme for Numerical Solution of Initial Value Problems in Ordinary Differential Equations. Pacific Journal of Science and Technology. 11(2): 260-267.

[Wds] D. Wood, D. Dimitrov and H. Kojouharov (2015). A nonstandard finite difference method for n-dimensional productive-destructive systems. *J. of Difference Equations and Applications*, 21(3): 240-254.

[WdsTh] Daniel Wood (2015). *Advancement and Applications of Nonstandard Finite Difference Methods*, Ph.D. Thesis, University of Texas at Arlington

[GnzP] Gilberto C. González-Parra Miguel Díaz-Rodriguez Victor Comezaquira (2012). A Nonstandard Finite Difference Scheme For An Epidemic Model Of Fractional Order. *Avances En Simulación Computacional Y Modelado Numérico* E. Dávila, G. Uzcátegui, M. Cerrolaza (Editores).

[D&D] Kushal Dhar Dwivedia and S. Das (2019). Numerical solution of the nonlinear diffusion equation by using non-standard/standard finite difference and Fibonacci collocation methods. *Eur. Phys. J. Plus* (2019) 134: 608.

[App] A. R. Appadu (2013). Numerical Solution of the 1D Advection-Diffusion Equation Using Standard and Nonstandard Finite Difference Schemes. *Journal of Applied Mathematics*, vol. 2013, Article ID 734374, 14 pages. https://doi.org/10.1155/2013/734374

[Swe] N. H. Sweilam, A. M. Nagy, L. E. Elfahri (2019). Nonstandard Finite Difference Scheme for the Fractional Order Salmonella Transmission Model

[Dav] E. Dávila, G. Uzcátegui, M. Cerrolaza (Editores, 2012). A Nonstandard Finite Difference Scheme For An Epidemic Model Of Fractional Order. *Avances En Simulación Computacional Y Modelado Numérico*

[B,CC]=Brauer, F. and Castillo-Chávez, C. (2001). *Mathematical Models in Population Biology and Epidemiology*. Springer

[Z&K] Zakharchenko, Vladimir D. and Kovalenko, Ilya G. (2018). Best approximation of the Fractional Derivative Operator by Exponential Series.
Reprinted from: *Mathematics* 2018, 6, 12, doi: 10.3390/math6010012

[Demp] Dempster, M.A.H. (2000). Wavelet Based PDE Valuation of Derivatives. 7[th] Annual CAM Workshop on Mathematical Finance, December 01, 2000

[B&P1] Bohner, M. and Peterson, A. (2001). *Dynamic Equations on Time Scales: An Introduction with Applications*. Birkhauser Boston Inc., Boston, MA, 2001

[Bhn] M. Bohner, S. Streipert and D.F.M. Torres (2019). Exact solution to a dynamic SIR model. *Nonlinear Analysis: Hybrid Systems* 32 228–238

[DPC-M 3] Clemence-Mkhope, D. P. (2021, Preprint). The Exact Spectral Derivative Discretization Finite Difference (ESDDFD) Method for the Fokker-Planck Equation

[CL&V] E. Capelas de Oliveira, S. Jaros, J. Vaz Jr. (2020). On the mistake in defining fractional derivative using a non-singular kernel. arXiv:1912.04422v3 [math.CA] 29 Jan 2020

[FM] [Mainardi, 2018] Mainardi, Francesco (2018). A Note on the Equivalence of Fractional Relaxation Equations to Differential Equations with Varying Coefficients.
Reprinted from: *Mathematics* 2018, 6,8, doi: 10.3390/math6010008

[Xin] Xin, Baogui, Wei Peng, Yekyung Kwon, and Yanqin Liu (2019). Modeling, discretization, and hyperchaos detection of conformable derivative approach to a financial system with market confidence and ethics risk. *Advances in Difference Equations* 2019:138

[Mhmmdz] Mohammadnezhad, Vahid, Mostafa Eslami, Hadi Rezazadeh (2020). Stability Analysis of Linear Conformable Fractional Differential Equations System with Time Delays; *Bol. Soc. Paran. Mat.* (3s.) v. 38 6: 159–171.

[K&G2] Kolwanker, K.M. and Gangal, A.D. (1997). Local fractional Fokker-Planck equation. *Physical Review Letters*. 80(2):214-217 (arXiv:cond-mat/9801138)

[O&M] Odibat, Zaid M and Momani, Shaher (2008). An Algorithm for the numerical solution of differential equations of fractional order. *J. Appl. Math. & Informatics*, No 1-2_15-27

[DPC-M 4] Clemence-Mkhope, D. P. (2021, Preprint). On the Conformable Exponential Function and the Conformable Euler's Finite Difference Method

[DPC-M 5] Clemence-Mkhope, D. P. (2021, Preprint). On the Generalized Euler's Finite Difference Method and the Conformal-Caputo Derivative Connection

[DPC-M 6] Clemence-Mkhope, D. P. (2021, Preprint). The Exact Spectral Derivative Discretization Finite Difference (ESDDFD) Method for Wave Models: Dimensional and Other Extensions

[Nkt] Robin Nittka (2010). *Elliptic and Parabolic Problems with Robin Boundary Conditions on Lipschitz Domains*. Dissertation zur Erlangung des Doktorgrades Dr. rer. nat. der Fakultät für Mathematik und Wirtschaftswissenschaften der Universität Ulm. Downloaded 04/29/2021 from: https://d-nb.info/1004249772/34[P&Z]

[K&T] *Katz, Mikhail; Tall, David (2011), Tension between Intuitive Infinitesimals and Formal Mathematical Analysis, Bharath Sriraman, Editor. Crossroads in the History of Mathematics and Mathematics Education.* The Montana Mathematics Enthusiast *Monographs in Mathematics* Education 12, Information Age Publishing, Inc., Charlotte, *NC,* arXiv*:1110.5747,* Bibcode*:2011arXiv1110.5747K*

[WT] Watkins, Thayer (undated). The Calculus of Infinitesimals. Downloaded 05/27 from https://www.sjsu.edu/faculty/watkins/infincalc.htm

[Tao] Tao, Terence (2020). "The Notorious Collatz conjecture". Retrieved from https://terrytao.files.wordpress.com/2020/02/collatz.pdf

[Rd20a] Reid, Fabian (2020). On the Pattern in the Collatz Problem. @inproceedings{Reid2020OnTP

[Rd20b] Reid, Fabian (2020; Draft Manuscript). The Visual Pattern and Proof of the Collatz Conjecture

[Rd21] Reid, Fabian (2021). The Visual Pattern in the Collatz Conjecture and Proof of No Non-Trivial Cycles. https://arxiv.org/abs/2105.07955

[MD] Mauldin, R. Daniel (1997). "A Generalization of Fermat's Last Theorem: The Beal Conjecture and Prize Problem" (PDF). *Notices of the AMS*. 44 (11): 1436–1439.